\documentclass[amstex,12pt,amssymb]{article}
\usepackage{mathtext}
\usepackage{amsmath}
\usepackage{amssymb}
\usepackage{amsxtra}
\usepackage{latexsym}
\usepackage{ifthen}

\usepackage{color}

\newcommand {\beq}{\begin{equation}}
\newcommand {\eeq}{\end{equation}}
\newcommand{\reff}[1]{(\ref{#1})}

\newtheorem{theorem}{Theorem}[section]
\newtheorem{lemma}[theorem]{Lemma}

\newtheorem{definition}[theorem]{Definition}
\newtheorem{remark}[theorem]{Remark}

\def\lam{\lambda} \def\veps{\varepsilon}
\def\diy{\displaystyle}

\begin{document}
\def\diy{\displaystyle} \def\vphi{\varphi}

\title{A remark on normalizations in a local principle of large deviations}

\author{ A.V.Logachov$^{\dagger ,\ast,\Diamond,\flat}$,
Y.M.Suhov$^{\sharp ,\ddagger ,\triangle}$, N.D.Vvedenskaya$^\sharp$, A.A.Yambartsev$^{\flat}$ }

\maketitle

 {\footnotesize

\noindent $^\sharp$ Dobrushin Laboratory, Institute for Information Transmission Problems, Russian
Academy of Science, Bolshoj Karetnyj Per 19, Moscow 127051, RF; Vvedenskaya's E-mail:\ \ ndv@iitp.ru

\noindent $^\dagger$ Laboratory of Probability Theory and Mathematical Statistics,
Sobolev Institute of Mathematics, Siberian Branch of the Russian Academy of Sciences,
Koptuga, 4, Novosibirsk 630090, RF; Logachov's E-mail:\ \ omboldovskaya@mail.ru

\noindent $^\ast$ Department of High Mathematics, Siberian State University of Geosystems and Technologies,  Ul Plahotnogo 10, Novosibirsk 630108, RF

\noindent $^\Diamond$  Novosibirsk State University of Economics and Management, Kamenskaya Ul 56,  Novosibirsk 630099, RF

\noindent $^\ddagger$  Math Department, Penn State University, McAllister Buid, University Park, State
College, PA 16802, USA

\noindent $^\triangle$ Statistical Laboratory, DPMMS, University of Cambridge, Wilberforce Rd,
Cambridge CB3 0WB; Suhov's E-mail:\ \ yms@statslab.cam.ac.uk

\noindent $^\flat$ Department of Statistics, Institute of Mathematics
and Statistics, University of S\~ao Paulo, Rua do Mat\~ao 1010, CEP 05508--090, S\~ao Paulo SP, Brazil;
Yambartsev's E-mail:\ \ yambar@ime.usp.br
}


\vspace{1cm}

\begin{abstract}
This work is a continuation of \cite{VLSY1}.
We consider a continuous-time birth-and-death process in which the transition rates
have an asymptotical power-law dependence upon the position of the process.
We establish rough exponential asymptotic for the probability that a sample path
of a normalized process lies in a neighborhood of a given nonnegative continuous
function. We propose a variety of normalization schemes for which the large deviation
functional preserves its natural integral form.
\end{abstract}



\section{Introduction}

The study of birth-and-death processes provides an interesting topic, both theoretically
and in a number of applications. As examples, we quote the information theory (encoding
and storage of information, see \cite{SS1}), biology and chemistry (models of growth and
extinction in systems with multiple components,  see \cite{SS2},  \cite{MSSZ}), and economics
(models of competitive production and pricing, \cite{MPY}, \cite{VSB}).

We consider a continuous-time Markov process $\xi (t)$, $t\geq 0$, with state space
$\mathbb{Z}^+:=\{0\}\cup\mathbb{N}$, and with $\xi(0) =0$.
The evolution of the process $\xi (\,\cdot\,)$ is governed by the transition rates
$\lambda(x)>0$ for the jump $x\to x+1$, $x\in\mathbb{Z}^+$, and
$\mu(x)>0$ for the jump $x\to x-1$, $x\in\mathbb{N}$. For $x=0$ we set $\mu (x)=0$.
We will work with events
that exclude an explosion of the process in a given time-slot $0\leq t\leq T$.

A key assumption is that
\beq \label{condtns1}
\lim_{x\to
\infty}\frac{\lambda(x)}{Y(x)}= \lim_{x\to\infty}\frac{\mu(x)}{Z(x)}=1,
\eeq
here
\beq \label{condtns2}
Y(x):=y(x)x^l,\;\;Z(x):=z(x)x^m,
\eeq
where  $l,m\geq 0, l\neq m$ (and
hence  $\max (l, m)>0$),  and  $y(x)$, $z(x)$ are the
\textit{slowly varying} functions at infinity.
(A function $a(x)$ is called slowly varying at infinity, if  for all $\beta >0$
$\lim\limits_{x\to\infty}\frac{a(\beta x)}{a(x)} =1$; see, for example \cite{Gal} for more details.)

We study properties of a normalized (scaled) process $\xi_{\vphi, T}$ where
\begin{equation}\label{xit}
\xi_{\vphi ,T}(t)=\frac{\xi(tT)}{\varphi(T)}, \ 0\leq t\leq 1.
\end{equation}
Here $T>0$ is parameter and $\vphi$ a positive function.
The conditions upon $\vphi$ is stated as follows:
\beq\label{condtns3}\begin{array}{l}
\lim\limits_{T\to\infty}\varphi (T)=\infty\;\hbox{ and }\;
\lim\limits_{T\rightarrow\infty}\diy\frac{\varphi(T)\ln \big(\varphi(T)\big)}{
TV(\varphi(T))}=0\\
\qquad\qquad\qquad\qquad\qquad
\hbox{ where }\;V(\varphi(T)):=\max\,\Big(Y(\varphi(T)),Z(\varphi(T))\Big). 
\end{array}\eeq
Under conditions \eqref{condtns1}, \eqref{condtns2}, \eqref{condtns3} we study the {\it local
large deviation principle}, i.e., the
logarithmic asymptotics for the
probability $\mathbf{P}(\xi_{\vphi ,T}\in U_\veps (f))$ that the path
of the scaled process
lies in an $\veps$-neighborhood $U_\veps (f)$ of a given continuous function $t\in [0,1]
\mapsto f(t)$ with $f(0)=0$ and $f(t)>0$ for $t>0$. More precisely,
we establish the existence of the limit
\beq\label{lldp}\begin{array}{l}
\lim\limits_{\varepsilon\rightarrow 0}\limsup\limits_{T\rightarrow \infty}\diy\frac{1}{\psi(T)}
\ln\mathbf{P}(\xi_{\vphi ,T}(\,\cdot\,)\in U_\varepsilon(f))\\
\qquad
=\lim\limits_{\varepsilon\rightarrow 0}\liminf\limits_{T\rightarrow \infty}\diy\frac{1}{\psi(T)}
\ln\mathbf{P}(\xi_{\vphi ,T}(\,\cdot\,)\in U_\varepsilon(f))=-I(f),\;\;\hbox{ where }
\psi (T) =TV(T).\end{array}\eeq
The point is that under the above formalism \eqref{condtns3}, \eqref{lldp}
the {\it large deviation functional} $I(f)$
does not depend on the choice of $\vphi$ and has a natural integral form:
\beq\label{rflldp}
I(f)=\int_0^1f^{l\vee m}(t)dt, 
\eeq
here and below $v\vee w$ stands for the maximum of positive numbers $v,w$. Next, for any
$f,g \in \mathbb{D}[0,1]$
\beq
\rho(f,g)=\sup\limits_{t\in[0,1]}|f(t)-g(t)|,
\eeq
and $\mathbb{D}[0,1]$ denotes the space of right-continuous functions with left-limit at each $t\in [0,1]$.

In an earlier paper by the authors \cite{VLSY1}, a similar result was proved for constant
functions $y(x)$, $z(x)$ and $\vphi (T)=T$. The present work is an attempt to answer
the question to what extent the result of  \cite{VLSY1} can be generalized without
changing the form of the functional $I(f)$. The second motivation comes from a comparison
with the case of constant values $\lam (x)=\lam $ and $\mu (x)=\mu$ (the latter for $x\geq 1$).
In our scheme, this happens when $l=m=0$.
Here, depending on the choice of the space-scaling factor $\vphi (T)$, one distinguishes
between moderate (when $\vphi (T)/{\sqrt T}\to\infty$ and $\vphi (T)/T\to 0$),
large (when $\vphi (T)/T\to C\in (0,\infty )$) and super-large 
(when $\vphi (T)/T\to\infty$) deviations, with
different forms of $I(f)$. It turns out that under the conditions introduced in the current
paper, the large deviation functional preserves its form regardless of the choice of
function $\vphi$.

The idea and the method of proof goes back to \cite{MPY, VLSY1, Log1}; this provides
certain limitations for the parameters of the scheme.
We would like to note that the case $l=m$ is not covered by our condition \eqref{condtns2}
and hence is not considered in this paper,
although it was included in \cite{VLSY1} in a more specific situation. \  (In some sense,  $l=m$
it is the most difficult case within the above formalism.) 

The paper is organized as follows: in Section 2 we introduce our main result (Theorem~\ref{th2.1})
and key lemmas: Lemma~\ref{l2.0} -- \ref{l2.4}. In Section 3 we prove Theorem~\ref{th2.1} and the lemmas.
In Section 4 we prove the auxiliary results. 

\section{Basic definitions and the main result}

We set
\beq\label{setF} F = \{ f\in\mathbb{C}[0,1]: \ f(0)=0 \text{ and } f(t)>0, 0<t\leq 1\}.\eeq

\begin{theorem} \label{th2.1}
Under conditions \eqref{condtns1}, \eqref{condtns2}, \eqref{condtns3}
the family of random processes $\{\xi_{\vphi,T}( t ),\;0\leq t\leq 1\}$ defined by \eqref{xit}
satisfies the {\rm{LLDP}} on the set $F$, with the normalized function $\psi (T)$ as in \eqref{lldp}
and the rate function as in \eqref{rflldp}:
$$\psi(T)=TV(\varphi(T)),\quad I(f)=\int_0^1f^{l\vee m}(t)dt.$$
\end{theorem}


Note that if $l\vee m>1$ and
$\lim\limits_{T\rightarrow\infty} \varphi(T) = \infty$ the condition \eqref{condtns3} obviously holds.

As in \cite{MPY, VLSY1}, we consider an auxiliary Markov process $\{\zeta(t),\;t\in[0,T]\}$, on
$\mathbb{Z}$, homogeneous in time and space
$\mathbb{Z}$, with rate $1$ and equiprobable ($1/2$) jumps $\pm 1$.
Denote by $X_T$ the set of piecewise-constant right-continuous
functions on the interval $[0,T]$ starting at zero with jumps $\pm 1$.

The first auxiliary statement is Lemma \ref{l2.0} below; we give it without proof
as it is straightforward.

\begin{lemma} \label{l2.0}
{\rm{(Cf. \cite{MPY,VLSY1}.)}} The distribution of the random process
$\xi(\,\cdot\,)$ on $X_T$ is absolutely continuous with respect to that of a process $\zeta(\,\cdot\,)$.
The corresponding density $\mathbf{p}=\mathbf{p}_T$ on $X_T$
(the Radon-Nikodym derivative) has the form:
\begin{equation}\label{8a}
\mathbf{p}(u)=\left\{ \begin{array}{ll}2^{N_T(u)}\Bigl(\prod\limits_{i=1}^{N_T(u)}
e^{-(h(u(t_{i-1}))-1)\tau_{i}}\nu(u(t_{i-1}),u(t_i))\Bigr)\\
\qquad \times e^{-(h(u(t_{N_T(u)})-1))(T-t_{N_T(u)})
},\;\;\qquad\quad\quad \;\; \mbox{ if } \;N_T(u)\geq 1,\\
e^{-(h(0)-1)T}, \qquad\qquad\qquad\qquad\qquad\qquad\qquad  \mbox{ if } \;  N_T(u)= 0.\\
\end{array} \right. 
\end{equation}
where $h(\cdot) := \lambda(\cdot) + \mu(\cdot)$. Here the function $u(\cdot)$  has $N_T(u)$ jumps at the moments   \
$t_1,t_2,...,t_{N_T(u)}$ \ such that \ $0=t_0<t_1<...<t_{N_T(u)}<T< t_{N_T(u)+1}$,
$\tau_i=t_i-t_{i-1}$. Further, $\nu (u(t_{i-1}),u (t_i))$ is given by
\beq\label{nu0}
\nu (u(t_{i-1}),u (t_i))= \left\{ \begin{array}{ll}\lambda(u (t_{i-1})), &
\mbox{ if }\; u (t_i)-u (t_{i-1})=1;\\
\mu(u (t_{i-1})), & \mbox{ if }\;  u (t_i)-u (t_{i-1})= -1. \end{array} \right.
\eeq
\end{lemma}

Let $N_T(\zeta)$ be the number of jumps of $\zeta(t)$ on the interval $[0,T]$.
The claim of Lemma~\ref{l2.0} is equivalent to the fact that for any measurable set $G\subseteq X_T$
\beq
\mathbf{P}(\xi(\cdot)\in G)=e^T \mathbf{E}\bigl( e^{-A_T(\zeta)}e^{B_T(\zeta)+N_T(\zeta)\ln2};
\zeta(\cdot)\in G \bigr). \label{3}\eeq
Here
\begin{eqnarray}
\label{4}
A_T(\zeta)&:=&\int_0^T h(\zeta(t))dt
\nonumber \\
&=&\left\{ \begin{array}{ll}\sum\limits_{i=1}^{N_T(\zeta)}h(\zeta(t_{i-1}))\tau_{i}+
h(\zeta(t_{N_T(\zeta)}))(T-t_{N_T(\zeta)}),&
\mbox{if }\; N_T(\zeta)\geq1;\\
h(0)T, & \mbox{if }\;  N_T(\zeta)=0. \end{array} \right.
\\
B_T(\zeta)&:=&\left\{ \begin{array}{ll}\sum\limits_{i=1}^{N_T(\zeta)}\ln(\nu(\zeta(t_{i-1}),\zeta(t_i))), &
\mbox{if }\; N_T(\zeta)\geq1;\\
0, & \mbox{if }\;  N_T(\zeta)=0. \end{array} \right. \label{5}
\end{eqnarray}
Below we use (\ref{3}) in the study of asymptotic behavior of
$\ln \mathbf{P}(\xi_{\vphi,T}(\cdot)\in U_\varepsilon(f))$.

The proof of Theorem~\ref{th2.1} shows that in the case $l\neq m$ the main contribution to
this asymptotic comes from $A_T(\zeta)$.

Consider the sequence of scaled processes
\beq\label{zetT}\zeta_{\vphi,T}(t)=\frac{\zeta(tT)}{\varphi(T)}, \ \ t\in[0,1].\eeq
Further on, we write, for brevity, $N_T, A_T, B_T$
 instead of $N_T(\zeta), A_T(\zeta), B_T(\zeta)$.

\smallskip

\begin{lemma}\label{l2.3} Let the conditions of Theorem {\rm{\ref{th2.1}}} be fulfilled. Then
$$
\lim\limits_{\varepsilon\rightarrow 0} \limsup\limits_{T\rightarrow\infty}
\frac{1}{TV(\varphi(T))}\ln\mathbf{E}\bigl(e^{B_T+N_T\ln2};\zeta_{\vphi,T}\in U_\varepsilon(f)\bigr)\leq0,
$$
where $f\in F$.
\end{lemma}

\begin{lemma} \label{l2.4} Under the conditions of Theorem~\ref{th2.1}, then
$$
\lim\limits_{\varepsilon\rightarrow 0} \liminf\limits_{T\rightarrow\infty} \frac{1}{TV(\varphi(T))}\ln\mathbf{E}\bigl(e^{B_T+N_T\ln2};\zeta_{\vphi,T}\in U_\varepsilon(f)\bigr)\geq0,
$$
where $f\in F$.
\end{lemma}

\section{Proof of Theorem~\ref{th2.1} and Lemmas~\ref{l2.3}, \ref{l2.4}}
\def\diy{\displaystyle}

\smallskip

{\sc Proof of Theorem~\ref{th2.1}.} First, let us estimate the term $A_T$
$$
A_T=\int_0^T h(\zeta(t))dt=T\int_0^1 h(\varphi(T)\zeta_{\vphi,T}(s))ds.
$$
 We consider a  set of trajectories $\zeta(\cdot)$ where
$\zeta_{\vphi,T}\in U_{\veps}(f)$.

For fixed
 $\veps$ let $\delta:=\delta(\veps)=\max\limits_{0\leq t\leq 1} \{t :f(t)<2\veps\}$.
 We note that $\lim\limits_{\veps \to 0}\delta=0$ for all functions from the set $F$.
By (\ref{condtns1}),\ (\ref{condtns2}) for any $\gamma_0>0$, $s\in [\delta,1]$ and
sufficiently large $T>0$
\beq\label{12}
1-\gamma_0 \leq \frac{h(\varphi(T)\zeta_{\vphi,T}(s))}{V(\varphi(T)) (\zeta_{\vphi,T}(s))^{l\vee m} }
\leq 1+\gamma_0.
\eeq
By (\ref{12}) for all sufficiently large $T$
\beq\label{15}
\begin{aligned}
T\int_\delta^1&(1-\gamma_0)V(\varphi(T))(f(s)-\varepsilon)^{l\vee m}ds  \leq \ A_T
\\
&  \leq T\int_0^\delta h(\varphi(T)\zeta_{\vphi,T}(s))ds
+T\int_\delta^1(1+\gamma_0)
V(\varphi(T))(f(s)+\varepsilon)^{l\vee m}ds.
\end{aligned}
\eeq
Thus, we get that
\beq\label{16}
\begin{aligned}
TV(\varphi(T)) & \int_\delta^1(1-\gamma_0)(f(s)-\varepsilon)^{l\vee m}ds
\ \leq \ A_T \\
& \leq T\delta (1+\gamma_0)V(3\varepsilon \varphi(T))+TV(\varphi(T))\int_\delta^1(1+\gamma_0)
(f(s)+\varepsilon)^{l\vee m}ds.
\end{aligned}
\eeq
Using   (\ref{3}) and the inequalities (\ref{16}), we shift to logarithms  obtaining  that

\beq\label{18}
\begin{aligned}
 -  
 \int_\delta^1 (1-\gamma_0)(f(s)-\varepsilon)^{l\vee m}ds+
\limsup\limits_{T\rightarrow\infty}\diy\frac{1}{TV(\varphi(T))}
\ln\mathbf{E} \bigl( e^{B_T+N_T\ln2};\zeta_{\vphi,T}\in U_{\varepsilon}(f) \bigr)\\
\geq\limsup\limits_{T\rightarrow\infty}\diy\frac{1}{TV(\varphi(T))}\ln\mathbf{P}(\xi_{\vphi,T}\in U_\varepsilon(f))\geq\liminf\limits_{T\rightarrow\infty}\displaystyle\frac{1}{TV(\varphi(T))}\ln\mathbf{P}
(\xi_{\vphi,T}\in U_\varepsilon(f))\\
\geq
-\int_\delta^1(1+\gamma_0)(f(s)+\varepsilon)^{l\vee m}ds
 +\liminf\limits_{T\rightarrow\infty}\diy\frac{1}{TV(\varphi(T))}
\ln\mathbf{E} \bigl( e^{B_T+N_T\ln2};\zeta_{\vphi,T}\in U_\varepsilon(f) \bigr)
\\
\phantom{\geq }
\ - \delta(1+\gamma_0)\frac{V(3\varepsilon \vphi(T))}{V(\varphi(T))}.
\end{aligned}
\eeq
Since (\ref{18}) is fulfilled for
any $\gamma_0>0$
as $\veps\rightarrow 0$, $\gamma_0\rightarrow0$ we receive
\beq\label{20}
\begin{aligned}
-\int_0^1f^{l\vee m}(s)ds+
\lim\limits_{\varepsilon\rightarrow 0}\limsup\limits_{T\rightarrow\infty}\diy\frac{1}{TV(\varphi(T))}\ln\mathbf{E}\bigl( e^{B_T+N_T\ln2};\zeta_{\vphi,T}\in U_\varepsilon(f) \bigr)\\
\qquad \geq\lim\limits_{\varepsilon\rightarrow 0}\limsup\limits_{T\rightarrow\infty}\diy\frac{1}{TV(\varphi(T))}\ln\mathbf{P}(\xi_{\vphi,T}(\cdot)\in U_\varepsilon(f))\\
\qquad \geq\lim\limits_{\varepsilon\rightarrow 0}\liminf\limits_{T\rightarrow\infty}\diy\frac{1}{TV(\varphi(T))}\ln\mathbf{P}(\xi_{\vphi,T}(\cdot)\in U_\varepsilon(f))\\
\qquad\geq-\int_0^1f^{l\vee m}(s)ds
+\lim\limits_{\varepsilon\rightarrow 0}\liminf\limits_{T\rightarrow\infty}\diy\frac{1}{TV(\varphi(T))}
\ln\mathbf{E} \bigr( e^{B_T+N_T\ln2};\zeta_{\vphi,T}\in U_\varepsilon(f) \bigr).
\end{aligned}
\eeq

Applying Lemmas~\ref{l2.3} and Lemma~\ref{l2.4} to inequalities (\ref{20})
finishes the proof of the theorem. $\Box$

\begin{remark} \label{r3.1} For the Yule pure birth process ($l>0$,  $\mu(x)\equiv 0$;
see for example {\rm{\cite{Fell}}}) the rate function has the form
$$
I(f)=\int_0^1f^l(t)dt, \ \ \ f\in F_M.
$$
Here $F_M$ is the set of continuous monotone increasing functions on $[0,1]$ starting from 0.
\end{remark}

\bigskip

{\sc Proof of Lemma~\ref{l2.3}.} In this lemma the goal is to establish the claimed
upper bound for the expected value
$\mathbf{E}\bigl(e^{B_T+N_T\ln2};\zeta_{\varphi ,T}\in U_\varepsilon(f)\bigr)$.  Obviously,
\beq\label{27}
\begin{aligned}
\mathbf{E} \bigl( e^{B_T+N_T\ln2};\zeta_{\varphi ,T}\in U_\varepsilon(f) \bigr) &:= E_{1}+E_{2},\;\hbox{ with }\\
\quad\quad  E_1&:=\mathbf{E} \bigl(e^{B_T+N_T\ln2};\zeta_{\varphi ,T}\in U_\varepsilon(f);
N_T\leq \Theta(T)  \bigr),\\
\quad\quad  E_2&:=\mathbf{E} \bigl( e^{B_T+N_T\ln2};\zeta_{\varphi ,T}\in U_\varepsilon(f);N_T >
 \Theta(T) \bigr),
\end{aligned}
\eeq
where $\Theta(T) := \sqrt{\frac{TV(\varphi(T))\varphi(T)}{\ln(\varphi(T))}}$.
%
Denote $M=\max\limits_{t\in[0,1]}f(t)\vee 1$. If $\zeta_{\varphi ,T}\in U_\varepsilon(f)$ and $N_T\leq  \Theta(T) $ then
for any $\gamma_1>0$ and for all sufficiently large $T$
$$
\begin{aligned}
B_T&=\sum\limits_{i=1}^{N_T}\ln \bigl( \nu(\zeta(t_{i-1}),\zeta(t_i)) \bigr) \\
&\leq  \Theta(T)
\Bigl( \ln \bigl(Y(\varphi(T))(M+\varepsilon)^l(1+\gamma_1) \bigr)+\ln \bigl(Z(\varphi(T))(M+\varepsilon)^m(1+\gamma_1)\bigr) \Bigr).
\end{aligned}
$$
Denote $k_1:=(M+\varepsilon)^{l+m}(1+\gamma_1)^2$.
As $V(\varphi(T))\leq y(\varphi(T))\vee z(\varphi(T))\varphi^{l\vee m}(T)$ and for sufficiently large $T$
$y(\varphi(T))\vee z(\varphi(T))\leq \varphi^{l\vee m}(T)$
we obtain the  inequality
\beq\label{E1}
\begin{aligned}
E_1 &\leq \exp\Bigl\{  \Theta(T)
\ln \bigl( k_1Y(\varphi(T))Z(\varphi(T)) \bigr) \Bigr\} 2^{ \Theta(T) }  \\
&\leq  \exp\Bigl\{ \Theta(T) \ln \bigl(2 k_1Y(\varphi(T))Z(\varphi(T)) \bigr) \Bigr\}\leq
\exp\Bigl\{ \Theta(T) \ln \bigl(2 k_1V^2(\varphi(T))\bigr) \Bigr\}
\\ &\leq\exp\Bigl\{ \Theta(T) \ln \bigl(2 k_1\varphi^{2(l\vee m)}(T)\bigr) \Bigr\}.
\end{aligned}
\eeq

\medskip

Next, denote by $k_+$ and $k_-$ the number of positive and negative jumps of the process $\zeta_{\varphi ,T}(\cdot)$ and let $L=k_+-k_-$.
For $\zeta_{\varphi ,T}(\cdot)\in U_\varepsilon(f)$ the following inequality holds
\beq
f(1)-\varepsilon \leq \zeta_{\varphi ,T}(1) \leq f(1)+\varepsilon.  \label{8}
\eeq
Since the jumps of the process $\zeta_{\varphi ,T}(\cdot)$ are  $\pm 1/\varphi(T)$,
by inequality (\ref{8}) we have
\beq
(f(1)-\varepsilon)\varphi(T) \leq L \leq (f(1)+\varepsilon)\varphi(T),  \label{9}
\eeq
and
\beq
k_+ + k_- =N_T, \ \ \  k_+=\frac{N_T+L}{2}, \ \ \ k_-=\frac{N_T-L}{2}. \label{10}
\eeq
As $\zeta_{\varphi ,T}\in U_\varepsilon(f)$, we obtain from (\ref{10}) that
for any $\gamma_1>0$ and for $T$  large enough,
\beq
\begin{aligned}
B_T=\sum\limits_{i=1}^{N_T}\ln \bigl( \nu(\zeta(t_{i-1}),\zeta(t_i)) \bigr) 
\leq \frac{ N_T+L}{2}\ln \bigl(Y(\varphi(T))(M+\varepsilon)^l(1+\gamma_1)\bigr)&
\\
+\frac{ N_T-L}{2}\ln\bigl(Z(\varphi(T))(M+\varepsilon)^m(1+\gamma_1)\bigr)&.
\end{aligned} \label{art1}
\eeq
Since
$N_T> 
 \Theta(T) $, we get, by using (\ref{9}) and the condition (\ref{condtns3}), that
\beq \label{10.11}
\lim_{T\to \infty} \frac{N_T}{L} = \infty.
\eeq
Thus, by (\ref{art1}) and (\ref{10.11}), for any $\gamma_1>0$ and all sufficiently $T$ we obtain
$$
\begin{aligned}
B_T & \leq \frac{N_T}{2} \ln \bigl(k_1Y(\varphi(T))Z(\varphi(T)) \bigr) + \frac{L}{2}\ln \Biggl( \frac{Y(\varphi(T))}{Z(\varphi(T))}
(M+\varepsilon)^{l-m} \Biggr) \\
&\leq \frac{N_T}{2}(1+\gamma_1) \ln\bigl(k_1Y(\varphi(T))Z(\varphi(T))\bigr).
\end{aligned}
$$
Hence,
\beq\begin{array}{rcl}
E_2 &\leq& \mathbf{E} \Bigl(e^{B_T+N_T\ln2};N_T\geq \Theta(T) +1\Bigr)\\
\qquad\qquad &\leq& \mathbf{E}\exp\Bigl\{\displaystyle\frac{N_T}{2}(1+\gamma_1)
\ln\bigl(4k_1Y(\varphi(T))Z(\varphi(T))\bigr)\Bigr\}. \label{art2}
\end{array}
\eeq
Since $ N_T$ has the Poisson distribution with parameter $T$,
$$
\mathbf{E}e^{\theta N_T}=e^{T(e^\theta-1)}.
$$
Therefore, from (\ref{art2}) it follows that
\beq
E_2\leq \exp\Bigl\{ k_2T
\bigl( Y(\varphi(T))Z(\varphi(T)) \bigr)^{(1+\gamma_1)/2}\Bigr\},  \label{art3}
\eeq
where $k_2:=(4k_1)^{(1+\gamma_1)/2}$.

Now let us choose $\displaystyle\gamma_1<\frac{|l-m|}{l+m}$. Using inequalities
(\ref{E1}), (\ref{art3}), condition \eqref{condtns3} and an obvious
inequality $\ln (c+d)\leq \ln (2(c\vee d))$, we obtain
$$
\begin{aligned}
&\lim\limits_{\varepsilon\rightarrow 0} \limsup\limits_{T\rightarrow\infty}
\frac{1}{TV(\varphi(T))}\ln\mathbf{E} \bigl(e^{B_T(\zeta)+N_T\ln2};\zeta_{\varphi ,T}\in U_\varepsilon(f) \bigr)
\\
&\leq\lim\limits_{\varepsilon\rightarrow 0} \limsup\limits_{T\rightarrow\infty}
\frac{ \Big[2\Theta(T)
\ln\bigl( k_1V^2(\varphi(T)) \bigr)\Big] \vee \Big[k_2T
\bigl( Y(\varphi(T))Z(\varphi(T)) \bigr)^{(1+\gamma_1)/2}\Big]}{TV(\varphi(T))}
\\
&\leq\lim\limits_{\varepsilon\rightarrow 0} \limsup\limits_{T\rightarrow\infty}
\left(\frac{ 2\sqrt{\varphi(T)}
\ln\bigl( k_1\varphi^{2(l\vee m)}(T) \bigr)}{\sqrt{TV(\varphi(T))\ln(\varphi(T))}}\vee
\frac{  k_2
\bigl( y(\varphi(T))z(\varphi(T))\varphi^{l+m}(T)\bigr)^{(1+\gamma_1)/2}}{V(\varphi(T))}\right)=0.\;\;\Box
\end{aligned}
$$

\bigskip

{\sc Proof of Lemma~\ref{l2.4}.}  The aim is to
lower-bound the term $\mathbf{E}\bigl(e^{B_T+N_T\ln2};\zeta_{\varphi ,T}\in U_\varepsilon(f)\bigr)$.
Set $k_3:=\inf\limits_{x\in \mathbb{Z}^+}\lambda(x)>0,$ and 
$k_4:=\inf\limits_{x\in \mathbb{N}}\mu(x)>0$. We note that $k_3>0$ and $k_4>0$.

Observe that if $\zeta_{\varphi ,T}\in U_\varepsilon(f)$ then $B_T\geq N_T\ln(k_3\wedge k_4)$,
where $v\wedge w$ is a minimum of positive numbers $v,w$.
Thus
\beq
\mathbf{E}\bigl(e^{B_T+N_T\ln2};\zeta_{\varphi ,T}\in U_\varepsilon(f)\bigr)\geq
\mathbf{E}\bigl(e^{N_T\ln(k_3\wedge k_4)};\zeta_{\varphi ,T}\in U_\varepsilon(f);N_T\leq C\varphi(T)\bigr), \label{26}
\eeq
where the constant $C>0$ depends on the function $f$ (see Lemma \ref{l4.1}) from the appendix.


$$
e^{C\varphi(T)\ln(k_3\wedge k_4)}\mathbf{P}(\zeta_{\varphi ,T}\in U_\varepsilon(f);N_T\leq C\varphi(T)).
$$

Thus, from (\ref{26})  it follows that
$$
\begin{aligned}
 \liminf\limits_{T\rightarrow\infty} &
 \frac{1}{TV(\varphi(T))}\ln\mathbf{E} \bigl( e^{N_T\ln(k_3\wedge k_4)};\zeta_{\varphi ,T}\in U_\varepsilon(f);N_T\leq C\varphi(T) \bigr)
\\
&\geq\liminf\limits_{T\rightarrow\infty}
 \frac{1}{TV(\varphi(T))}\ln\mathbf{P} \bigl( \zeta_{\varphi ,T}\in U_\varepsilon(f);N_T\leq C\varphi(T) \bigr).
\end{aligned}
$$

By Lemma~\ref{l4.1} from the appendix, 
$$
\liminf\limits_{T\rightarrow\infty}
 \frac{1}{TV(\varphi(T))}\ln\mathbf{P}\bigl(\zeta_{\varphi ,T}\in U_\varepsilon(f);N_T\leq C\varphi(T)\bigr)=0.
$$
This completes the proof of Lemma~\ref{l2.4}. $\Box$

\bigskip

\section{Appendix}


\begin{lemma} \label{l4.1} Let the condition \eqref{condtns1} be fulfilled. Then for any function
$f\in F$ there exists a constant $C\ (C=C(\varepsilon))$ such that for any $\varepsilon>0$
$$
\liminf\limits_{T\rightarrow\infty}
 \frac{1}{TV(\varphi(T))}\ln\mathbf{P} \bigl(\zeta_{\varphi ,T}\in U_\varepsilon(f);N_T\leq C\varphi(T)\bigr)=0.
$$
\end{lemma}

\smallskip

{\sc Proof.} The process $\zeta(t)$ can be represented as
$$
\zeta(t)=\zeta^{(1)}(t)-\zeta^{(2)}(t),
$$
where $\zeta^{(1)}(t)$ and $\zeta^{(2)}(t)$ are independent Poisson processes with parameter
$\mathbf{E}\zeta^{(1)}(t)=\mathbf{E}\zeta^{(2)}(t)=t/2$.

Since  $f$ is continuous there exists a continuous function of finite variation
$g$ such that $\rho(f,g)<\varepsilon/2$, $g(0)=0$. Moreover,  there exist
continuous monotone non-decreasing functions $g_+$ and $g_-$ such that
$$
g(t)=g_+(t)-g_-(t), \ \ \ g_+(0)=g_-(0)=0.
$$
Because of independence of  processes $\zeta^{(1)}$ and $\zeta^{(2)}$ we can write
\begin{equation}
\begin{aligned}
\mathbf{P}\bigl( \zeta_{\varphi ,T}(\cdot)\in U_\varepsilon(f); & N_T\leq C\varphi(T) \bigr) \\
{}\geq\mathbf{P} \bigl(\zeta^{(1)}_{\varphi ,T}(\cdot)&\in U_{\varepsilon/4}(g_+);N^{(1)}_T
\leq C_1\varphi(T)\bigr)\\
&{} \times
\mathbf{P} \bigl(\zeta^{(2)}_{\varphi ,T}(\cdot)\in U_{\varepsilon/4}(g_-);N^{(2)}_T
\leq C_2\varphi(T)\bigr) \ =:\ P_1P_2.\end{aligned}\end{equation}
Here, in analogy to \eqref{zetT},
$$\zeta^{(1)}_{\varphi ,T}(t)=\frac{\zeta^{(1)}(tT)}{\varphi(T)},\quad \zeta^{(2)}_{\varphi ,T}
(t)=\frac{\zeta^{(2)}(tT)}{\varphi(T)}.$$
Furthermore, $N^{(i)}_T$ stands for the number of jumps in $\zeta^{(i)}$ on $[0,T]$, $i=1,2$.
Finally,
$$C_1=g_+(1), \ \ \ C_2=g_-(1), \ \ \ C=C_1+C_2.$$

To lower-bound the probability $P_1$,
consider a partition of the unit interval by points $0=t_0<t_1< ... <t_K=1$  such that
$$\max\limits_{i=1,\dots,K}
(g_+(t_i)-g_+(t_{i-1}))<\frac{\varepsilon}{8}.$$
Since $\zeta^{(1)}$ is a process with independent increments, we get that for a sufficiently
large $T$
$$\begin{aligned}
P_1&\geq \prod\limits_{i=1}^K
\mathbf{P} \bigl( \zeta^{(1)}(Tt_i)-\zeta^{(1)}(Tt_{i-1})= \lfloor (g_+(t_i)-g_+(t_{i-1}))\varphi(T) \rfloor \bigr)
\\
&=\prod\limits_{i=1}^K \frac{e^{-T(t_i-t_{i-1})/2}(T(t_i-t_{i-1})/2)^{ \lfloor (g_+(t_i)-g_+(t_{i-1}))\varphi(T) \rfloor }}{ \lfloor (g_+(t_i)-g_+(t_{i-1}))\varphi(T) \rfloor !}
\\
&\geq \prod\limits_{i=1}^K \exp\left\{- \frac{T(t_i-t_{i-1})}{2} -
(g_+(t_i)-g_+(t_{i-1}))\varphi(T)\ln \bigl( (g_+(t_i)-g_+(t_{i-1}))\varphi(T) \bigr) \right\}
\\
&\geq \prod\limits_{i=1}^K \exp\left\{- \frac{T(t_i-t_{i-1})}{2}-
(g_+(t_i)-g_+(t_{i-1}))\varphi(T)\ln\bigl( g_+(1)\varphi(T) \bigr) \right\}
\\
&\geq \exp\left\{-T-
g_+(1)\varphi(T)\ln(g_+(1)\varphi(T))\right\},
\end{aligned}
$$
where $\lfloor b \rfloor$ is the integer part of the number $b$.

In the same way we obtain a lower bound for $P_2$:
$$P_2\geq \exp\left\{-T-
g_-(1)\varphi(T)\ln\bigl( g_-(1)\varphi(T) \bigr)\right\}.$$
Then from (\ref{condtns3}) it follows that
$$
\begin{aligned}
&\liminf\limits_{T\rightarrow\infty}
\ln\mathbf{P}\bigl( \zeta_{\varphi ,T}\in U_\varepsilon(f);N_T\leq C\varphi(T) \bigr)\geq
\liminf\limits_{T\rightarrow\infty}
\ln(\mathbf{P}_1\mathbf{P}_2)\quad
\\
&\quad\geq\liminf\limits_{T\rightarrow\infty}
 \frac{-2T-
(g_-(1)+g_+(1))\varphi(T)\ln\bigl((g_-(1)+g_+(1))\varphi(T)\bigr)}{TV(\varphi(T))}=0.
\end{aligned}
$$
This completes the proof of Lemma \ref{l4.1}.$\Box$

\section*{Acknowledgments}

We woulds like to
stress the role of E.A.Pecherskiy in the formulation of the problems leading to the current results.

NDV thanks Russian Science Foundation for the financial support through Grant
14-50-00150. AVL thanks FAPESP (S\~{a}o Paulo Research Foundation) for the financial
support via Grant 2017/20482
and also thanks RFBR (Russian Foundation for Basic Research) grant 18-01-00101.
YMS thanks The Math Department, Penn State University, for hospitality and support
and St John's College, Cambridge, for financial support.
AAY thanks CNPq (National Council for Scientific and Technological Development)
and  FAPESP for the financial support via
Grants 301050/2016-3 and 2017/10555-0, respectively.

\end{document}